\documentclass[12pt]{amsart}
\usepackage{latexsym,amssymb}
\usepackage{graphicx}
\usepackage{amsthm}
\usepackage{mathrsfs}
\usepackage{amsmath}
\usepackage{caption}

\numberwithin{equation}{subsection}

\newcommand{\tends}[1]{\mbox{\space \raise-2mm
\hbox{$\textstyle\longrightarrow\atop\scriptstyle {#1}$} \space}}

\newcommand{\vt}{\vec}

\date{\today}

\newtheorem{thm}{Theorem}[section]
 \newtheorem*{reftheorem1}{Theorem \reftotheorem}

 \theoremstyle{definition}
 
 \theoremstyle{remark}

 \newtheorem{notat}[thm]{Notation}

\let\oldmarginpar\marginpar
  \renewcommand\marginpar[1]
  {\-\oldmarginpar[\raggedleft\footnotesize\ \textit{#1}]%
  {\raggedright\footnotesize \textit{#1}}}
\begin{document}

\title[The projection of the vertices of an equilateral triangle]
{On the projection of the three vertices of an equilateral triangle}
\author{Maria Gabriella Kuhn}
\address{M. Gabriella Kuhn \\ Dipartimento di Matematica e Applicazioni \\ Universit\'a di Milano Bicocca\\ via Cozzi 55 \\ 20125 Milano\\ Italy}
\email{mariagabriella.kuhn@unimib.it}
\author{Silvio Norberto Riccobon}
\address{Scientific Committee of  \\ Matematica Senza Frontiere }
\email{sriccobon@gmail.com}

\thanks{}

\subjclass[2020]{Primary: 51M04; Secondary 12D10 }

\keywords{}

\date{}

\dedicatory{}

\begin{abstract}
In this paper we shall prove the following statement:
Given any three distinct points on a straight line r, there exist an equilateral triangle, whose circumcenter lies on r, such that the projections of its vertices on r are exactly the three given points.
\end{abstract}
\maketitle

\section{Introduction}
Consider the algebraic equation $x^3+px+q=0$.
It has been known since the 16th century from Cardano’s formulas that the equation admits three real roots if and only if $(\dfrac{q}{2})^2 + ( \dfrac{p}{3})^3 \leq 0$.
In this article we want to give a different approach to the study of the roots of this equation. That is, we will use a geometric theorem to produce formulas for the roots that do not involve complex numbers.
\section{The signed measurement of oriented segments}
\begin{flushleft}
In this paper we shall make intensive use of the following conventions:\\
$AB$ will denote the segment of endpoints $A$ and $B$.\\
$\vt{AB}$ will denote the segment oriented from $A$ to $B$.\\
Given an oriented line $r$ and two arbitrary points $A$ and $B$ on it, we define the signed measure of the oriented segment $\vt{AB}$ as the measure of the segment $AB$ preceded by the $+$ sign if the direction from $A$ to $B$ coincides with the direction of the straight line, with the  sign $-$ otherwise.\\
The measure of the segment $AB$ will be denoted $\overline{AB}$\\
The signed measure of the oriented segment $\vt{AB}$ will be denoted $(AB)$\\
Properties: 	The following properties hold:\\
\end{flushleft}
\begin{equation} \label{eq0_01}
(AB) = - (BA)
\end{equation} 
\begin{equation} \label{eq0_02}
(AB) = (AC) + (CB) \qquad \forall \quad C \in r
\end{equation} 

\maketitle

\section{Theorem of the three aligned points}
Statement: 	Given any three distinct points A, B, C, aligned on the straight line r, there exists an equilateral triangle with circumcenter at $O \quad  [ O \in r  \quad \land \quad (OA)+(OB)+(OC)=0 ]$ such that the projections of its vertices on r coincide with the given points $A, B, C.$\\
\begin{figure}[h]
\centering
\includegraphics[width=1\textwidth]{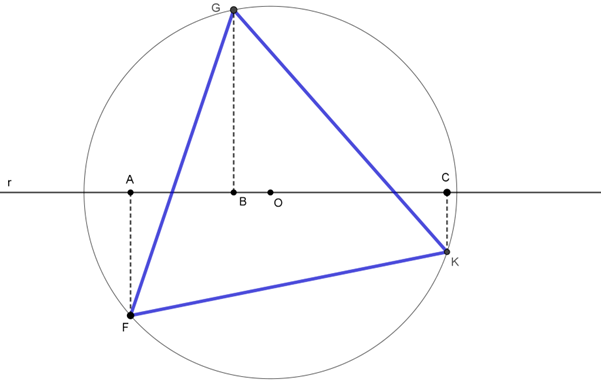}

\centering
\caption*{Figure 0: What we have to prove}

\end{figure}

N.B. 	The triangle that satisfies the theorem is not unique: its symmetric with respect to the straight line r also satisfies the theorem.
\\
\begin{figure}[h]
\centering
\includegraphics[width=0.6\textwidth]{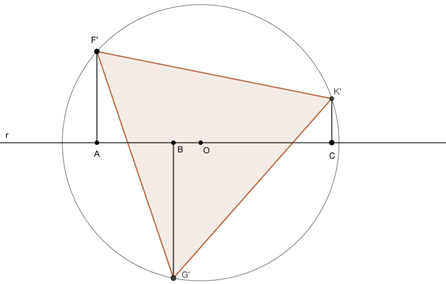}
\caption*{Figure 0.1 : The Symmetric Solution}
\end{figure}
\maketitle
\centering
PROOF\\
The proof is divided into 4 steps.
\centering
\subsection{1st Step}
\ \\
Geometric construction of the circumcenter O.\\
\begin{flushleft}
The first goal is to construct a point $O$ on the line r such that:\\
\begin{equation} \label{eq1_01}
(OA)+(OB)+(OC)=0
\end{equation} 
The condition \eqref{eq1_01} is equivalent to:\\
$(OC)=-(OA)-(OB)$\\
$(OC)=(AO)+(BO)$\\
$(OC)=(AC)+(CO)+(BC)+(CO)$\\
$(OC)=(AC)+(BC)+2(CO)$\\
$(OC)=(AC)+(BC)-2(OC)$\\
$3(OC)=(AC)+(BC)$\\
\  \\
\begin{equation} \label{eq1_02}
(OC)= \dfrac{(AC)+(BC)}{3}
\end{equation} 
\  \\
We proceed with the construction of a point O satisfying the condition \eqref{eq1_02} , which is equivalent to \eqref{eq1_01}.
Without loss of generality we may assume that the points A, B, C on r are in the following order:\\
\begin{figure}[h]
\centering
\includegraphics[width=1\textwidth]{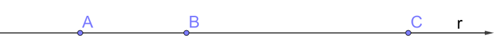}
\caption*{Figure 1.1}
\end{figure}
\  \\
Chose a point D on r so that (DA)=(BC) (see figure 1.2)
\begin{figure}[h]
\centering
\includegraphics[width=1\textwidth]{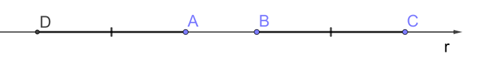}
\caption*{Figure 1.2}
\end{figure}
\  \\
By \eqref{eq0_02} we have:\\
\  \\
\begin{equation} \label{eq1_03}
(DC)=(DA)+(AC)=(BC)+(AC)
\end{equation} 

\begin{figure}[h]
\centering
\includegraphics[width=1\textwidth]{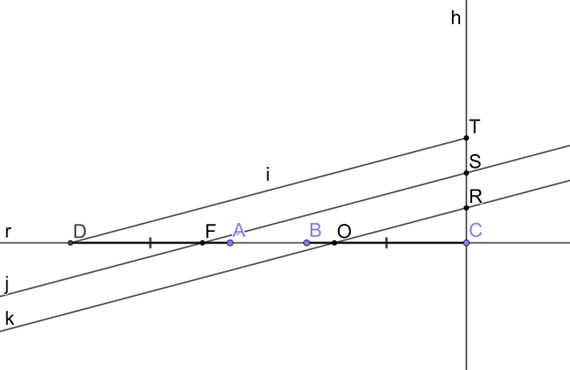}
\caption*{Figure 1.3}
\end{figure}
\  \\
From point $C$ draw the straight line $h$ perpendicular to $r$ (see figure 1.3).\\
Set a point $R$ on $h$ and construct the points $S, T$ so that\\
$\overline{CR}=\overline{RS}=\overline{ST} $.\\
Join $T$ with $D$ (segment $i$).\\
From $S$ draw the line $j$ parallel to the segment $i$, which intersects $r$ at the point $F$.\\
From $R$ draw the line $k$ parallel to the segment $i$, which intersects $r$ at the point $O$.\\
By Thales' theorem it results:\cite{ref1}
\ \\
$\overline{DF}=\overline{FO}=\overline{OC}$,\quad $\overline{DC}=\overline{DF}+\overline{FO}+\overline{OC}=3\overline{OC} $,\quad $\overline{OC}=\dfrac{\overline{DC}}{3}$\\
\  \\
Since $OC$ and $DC$ are both oriented from left to right, it results:\\
$(OC)=\dfrac{(DC)}{3}$\\
and by \eqref{eq1_03}: $(OC)=\dfrac{(AC)+(BC)}{3}$, as required by \eqref{eq1_02}\\
\end{flushleft}
\centering
\subsection{2nd Step}
\ \\
Based on the construction of figure 2, prove that $\overline{EO}=\overline{OC}$\\
\begin{figure}[h]
\centering
\includegraphics[width=1\textwidth]{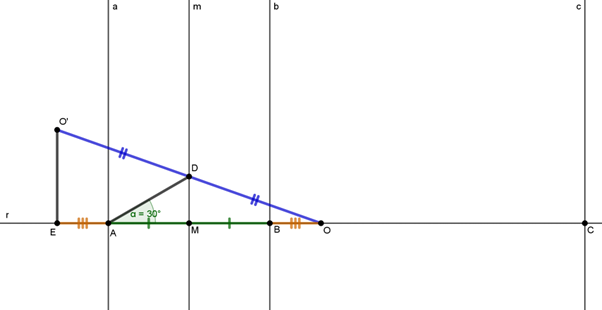}
\caption*{Figure 2}
\end{figure}
\begin{flushleft} 
We have constructed the point O so that:
\begin{equation} \label{eq2_01}
(OA)+(OB)+(OC)=0
\end{equation}
\  \\
In our case, with point $B$ closer to point $A$ than to $C$, is equivalent to (see \eqref {eq2_01}):
\begin{equation} \label{eq2_02}
\overline{OA}+\overline{OB}=\overline{OC}
\end{equation}
\  \\
We construct the midpoint $M$ of the segment $AB$, therefore
\begin{equation} \label{eq2_03}
\overline{AM}=\overline{MB}
\end{equation}
\  \\
From the points $A$, $M$, $B$, $C$ we draw the straight lines $a$, $m$, $b$, $c$ respectively, all perpendicular to the line $r$.\\ 
From point $A$ we draw a ray $s$ such that the angle between $r$ and $s$ is $\dfrac{\pi}{6}$. Let $D$ be intersection of $s$ and $m$.\\
\ \\
Construct the point O' symmetrical of O with respect to D, so that:
\begin{equation} \label{eq2_04}
\overline{OD}=\overline{DO'}
\end{equation}
\  \\
Let $E$ be the foot of the perpendicular drawn from $O'$ to the line $r$.\\
\ \\
By Thales' theorem we have that: $\dfrac{\overline{OM}}{\overline{ME}}=\dfrac{\overline{OD}} {\overline{DO'}}$
therefore $\overline{OM}=\overline{ME}$\\
\ \\
since then $\overline{EA}=\overline{EM}-\overline{AM} \quad \land \quad \overline{BO}=\overline{MO}-\overline{MB}$
differences of congruent segments, it results:\\
\begin{equation} \label{eq2_05}
\overline{EA}=\overline{BO}
\end{equation}
\  \\
We then have that $\overline{EO}=\overline{EA}+\overline{AO}=\overline{BO}+\overline{AO}$ therefore on the basis of \eqref {eq2_02}:\\
\begin{equation} \label{eq2_06}
\overline{EO}=\overline{OC}
\end{equation}
\end{flushleft}
\  \\
\centering
\subsection{3rd Step}
\ \\
Based on the construction of figure 3, prove that $\overline{OG}=\overline{OO'}$\\
\begin{figure}[h]
\centering
\includegraphics[width=1\textwidth]{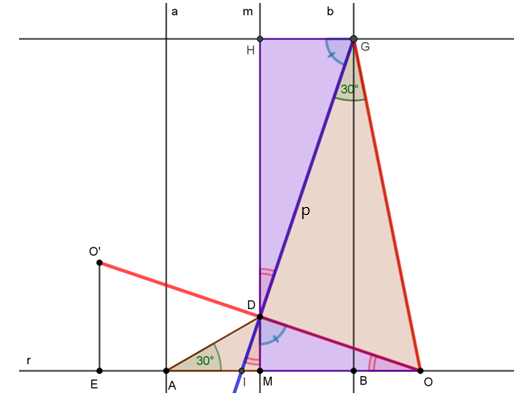}
\caption*{Figure 3}
\end{figure}
\begin{flushleft} 
Draw $p$, the axis of the segment $OO'$; let $G$ and $I$ be the intersections of $p$ with $b$ and with $r$ respectively.\\
$\measuredangle IDM \equiv \measuredangle GDH$ because angles opposite the vertex;\\
$\measuredangle BGD \equiv \measuredangle GDH$ because alternate internal angles between parallel lines; hence $\measuredangle IDM \equiv \measuredangle BGD$.\\
$\measuredangle MDO \equiv \measuredangle HGD$ because they are complementary of congruent angles.\\
Therefore the right-angled triangles $MOD$ and $HDG$ are similar and in particular there is a proportionality between the hypotenuses and the major cathetus:\\
\begin{equation} \label{eq3_01}
\dfrac{\overline{DO}}{\overline{DG}}=\dfrac{\overline{DM}}{\overline{HG}}
\end{equation}
\ \\
$\overline{HG}=\overline{MB}$ because opposite sides of a rectangle;\\
$\overline{AM}=\overline{MB}$ (see \eqref{eq2_03}), therefore $\overline{HG}=\overline{AM}$.\\
Based on this and on \eqref{eq3_01}we get: $\dfrac{\overline{DO}}{\overline{DG}}=\dfrac{\overline{DM}}{\overline{AM}}$.\\
\ \\
This demonstrates the similarity between right-angled triangles $DOG$ and $DMA$, whereby $\measuredangle DGO = \measuredangle MAD$; based on these facts, since $\measuredangle DGO=\dfrac{\pi}{6}$ , the hypotenuse OG measures twice the minor cathetus DO:\\
\ \\
$\overline{OG}=2\overline{DO}$.\\
\ \\
From \eqref{eq2_04} we get $\overline{OO'}=2\overline{DO}$, therefore:\\
\ \\
\begin{equation} \label{eq3_02}
\overline{OG}=\overline{OO'}
\end{equation}
\end{flushleft}
\ \\
\centering
\subsection{4th Step}
\ \\
\begin{figure}[h]
\centering
\includegraphics[width=1\textwidth]{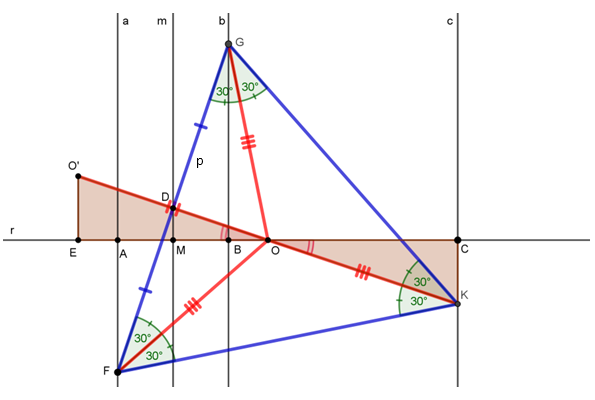}
\caption*{Figure 4}
\end{figure}
\begin{flushleft}
Consider the right-angled triangles $EOO'$ and $COK$:\\
$\overline{EO}=\overline{OC}$\quad see \eqref{eq2_06};\\
$\measuredangle EOO' \equiv \measuredangle COK$  because angles opposite the vertex.\\
Therefore, according to the criterion of congruence of right-angled triangles, triangles EOO' and COK are congruent.\\ Particularly:
\begin{equation} \label{eq4_01}
\overline{OK}=\overline{OO'}
\end{equation}\ \\
Let F be the intersection of $p$ with $a$.\\
Consider right-angled triangles DOF and DOG:\\
\ \\
by Thale’s theorem $\dfrac{\overline{DF}}{\overline{DG}}=\dfrac{\overline{MA}}{\overline{MB}}$, so that:
\begin{equation} \label{eq4_02}
\overline{DF}=\overline{DG}
\end{equation}\ \\
and the two right-angled triangles DOF and DOG are congruent. Particularly:
\begin{equation} \label{eq4_03}
\overline{OF}=\overline{OG}
\end{equation}
\begin{equation} \label{eq4_04}
\measuredangle DGO = \measuredangle DFO =\dfrac{\pi}{6}
\end{equation}\ \\
By the relations \eqref{eq3_02}, \eqref{eq4_01}, \eqref{eq4_03}, we have that:
\begin{equation} \label{eq4_05}
\overline{OF}=\overline{OG}=\overline{OK}
\end{equation}\ \\
The triangle FKG is isosceles \cite{ref2} since the perpendicular $KD$, conducted from the vertex $K$ to the base $FG$, divides it into two congruent parts \eqref{eq4_02}; therefore its height $KD$ is also the bisector, i.e
\begin{equation} \label{eq4_06}
\measuredangle OKF = \measuredangle OKG
\end{equation}\ \\
$\measuredangle OFK = \measuredangle OKF$ and $\measuredangle OGK= \measuredangle OKG$ because angles at the base of two isosceles triangles \label{eq4_05}.\\
By the transitive property of congruence, we have that:
\begin{equation} \label{eq4_07}
\measuredangle OFK = \measuredangle OKF = \measuredangle OKG = \measuredangle OGK
\end{equation}\ \\
\ \\
Since the sum of the internal angles of a triangle is $\pi$, referring to the triangle FGK, we have that:\\
\hspace{3cm}$\measuredangle FGK + \measuredangle GKF + \measuredangle KFG = \pi$\\
that is
\begin{equation} \label{eq4_08}
\measuredangle FGO+\measuredangle OGK+\measuredangle GKO+\measuredangle OKF+\measuredangle KFO+\measuredangle OFG=\pi
\end{equation}\ \\
\ \\
Based on \eqref{eq4_04} : 
\begin{equation*}
\measuredangle FGO +\measuredangle OFG=\dfrac{\pi}{3}
\end{equation*}
therefore \eqref{eq4_08} becomes:
\begin{equation*}
\measuredangle OGK+\measuredangle GKO+\measuredangle OKF+\measuredangle KFO=\dfrac{2\pi}{3}
\end{equation*}
\ \\
the four angles appearing on the first member are congruent \eqref{eq4_07}, therefore
\begin{equation*}
\measuredangle OGK=\measuredangle GKO=\measuredangle OKF=\measuredangle KFO=\dfrac{\pi}{6}
\end{equation*}
\ \\
from which we get:
\begin{equation*}
\measuredangle FGK=\measuredangle GKF=\measuredangle KFG=\dfrac{\pi}{3}
\end{equation*}
\ \\
Hence the triangle FGK is equilateral and its circumcenter coincides with O.\\
\centering
\ \\
q.e.d. 
\end{flushleft}
\section{Corollary}
\begin{flushleft}
It is known from Cardano’s formula that the equation
$x^3+px+q=0$ has three real roots if and only if: $p\leq -3\sqrt[3]{(q/2)^2}$.\\ 
Here we produce an alternative expression for its roots: namely\\
\ \\
$x_{1.2.3}=2\sqrt{-\dfrac{p}{3}} \cos [\dfrac{1}{3} \arccos(\dfrac{3q}{2p} \sqrt{-\dfrac{3}{p}})+\dfrac{2k\pi}{3}]$\quad with $k=-1,0,1$\\
\ \\
\ \\
PROOF\\
Every 3rd degree equation  $aX^3+ bX^2+cX+d=0$  , through the substitution:\\ 
\ \\
$X=x-\dfrac {b}{3a}$ , can be reduced to the form:
\begin{equation} \label{eq5_01}
x^3+px+q=0
\end{equation}\ \\
where  $p=\dfrac{3ac-b^2}{3a^2} \quad \land \quad q=\dfrac{2b^3-9abc+27a^2 d}{(27a^3}$\\
\ \\  
The coefficient of the 2nd degree term  is zero, since $x_1+ x_2+ x_3=0$ .\\
As we have already proved, the three real roots can be represented as the projections on the real line of the vertices of an equilateral triangle,
\begin{figure}[h]
\centering
\includegraphics[width=1\textwidth]{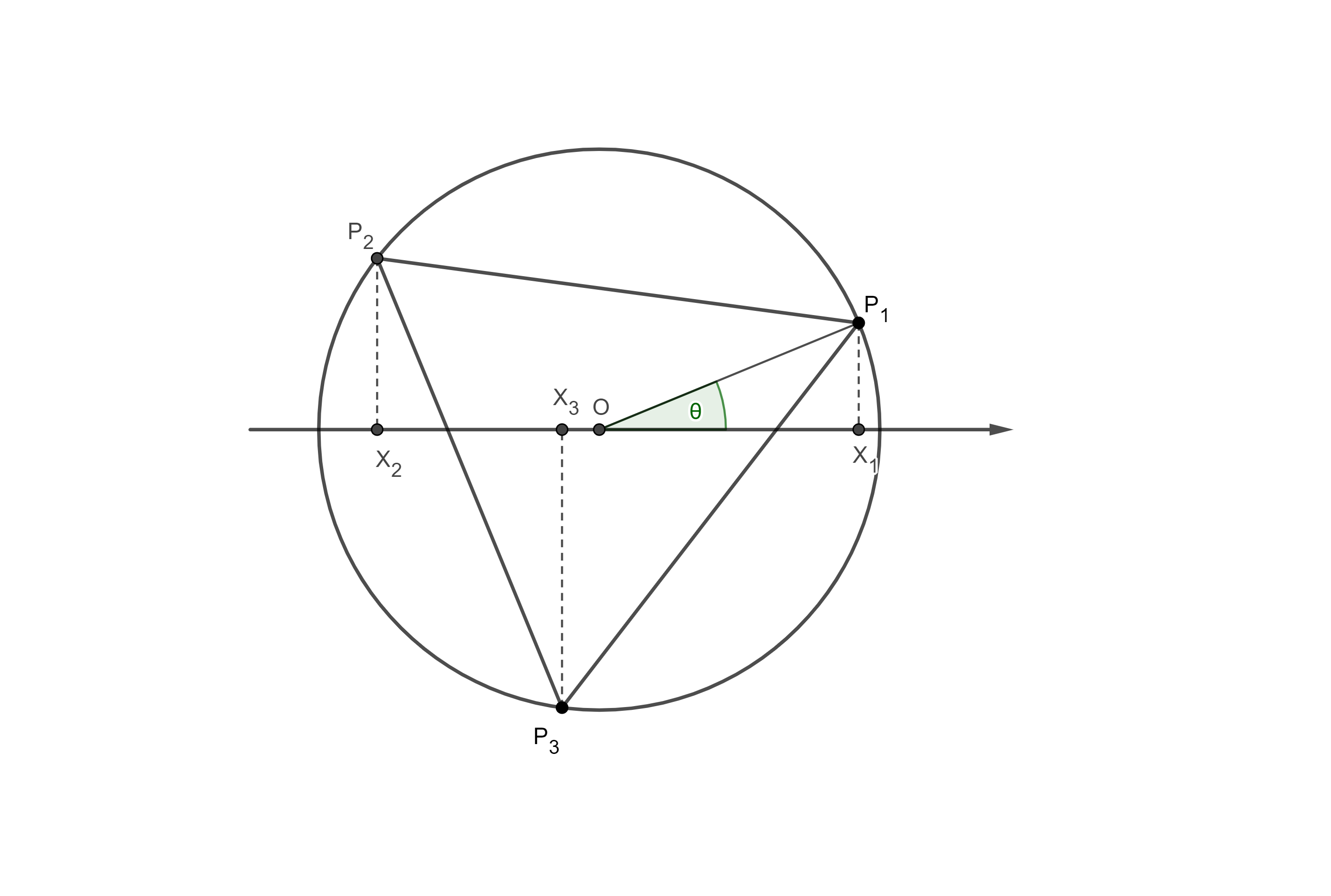}
\caption*{Figure 5}
\end{figure}
hence we may assume that they have the following expression:
\begin{equation*}
x_1=R \cos \theta  ,\qquad x_2=R \cos{(\theta+\dfrac{2\pi}{3})}  ,\qquad x_3=R \cos{(\theta-\dfrac{2\pi}{3})}
\end{equation*}
\ \\
So that you can write equation \eqref{eq5_01} in the form:
\begin{equation} \label{eq5_02}
[x-R\cos \theta ][x-R \cos{(\theta+\dfrac{2\pi}{3})}][x-R \cos{(\theta-\dfrac{2\pi}{3})} ]=0
\end{equation}\ \\
\begin{equation*}
\Longrightarrow \quad [x-R\cos \theta ][x+\dfrac{1}{2}R \cos \theta+\dfrac{\sqrt{3}}{2}R \sin \theta][x+\dfrac{1}{2}R \cos \theta-\dfrac{\sqrt{3}}{2}R \sin \theta]=0
\end{equation*}\ \\
\begin{equation*}
\Longrightarrow \quad  (x-R\cos \theta )(x^2+Rx \cos \theta+\dfrac{1}{4}R^2 \cos^2 \theta-\dfrac{3}{4}R^2 \sin^2 \theta)=0
\end{equation*}\ \\
\begin{equation*}
\Longrightarrow \quad (x-R\cos \theta )(x^2+Rx \cos \theta+R^2 \cos^2 \theta+\dfrac{1}{4}R^2 \cos^2 \theta-\dfrac{3}{4}R^2+\dfrac{3}{4}R^2 \cos^2 \theta=0
\end{equation*}\ \\
\begin{equation} \label{eq5_03}
\quad  x^3-\dfrac{3}{4}R^2x-R^3 \cos^3 \theta+\dfrac{3}{4}R^3\cos \theta=0
\end{equation}\ \\
\ \\
Comparing \eqref{eq5_01} with \eqref{eq5_03} we get:\\
\begin{equation*}
\left\{
\begin{array}{l}
p=-\dfrac{3}{4}R^2\\
q=-R^3 \cos^3 \theta+\dfrac{3}{4}R^3\cos \theta
\end{array}
\right.
\end{equation*}
\ \\
$\Longrightarrow$
\begin{equation*}
\left\{
\begin{array}{l}
R=2 \sqrt{-\dfrac{p}{3}}\\
q=-\dfrac{R^3}{4} (4 \cos^3 \theta-3\cos \theta)
\end{array}
\right.
\end{equation*}
\ \\
\ \\
From the first equality $R=2 \sqrt{-\dfrac{p}{3}}$ you get the condition:\\
\begin{equation} \label{eq5_04}
p \leq 0
\end{equation}\ \\
\ \\
\begin{equation*}
\Longrightarrow \quad q=-2 \sqrt{-(\dfrac{p}{3}})^3 (4 \cos^3 \theta-3\cos \theta)
\end{equation*}\ \\
\begin{equation*}
\Longrightarrow \quad  q=\dfrac{2p}{3} \sqrt{-\dfrac{p}{3}} (4 \cos^3 \theta-3\cos \theta) \qquad since \quad \lvert p \rvert=-p
\end{equation*}\ \\
\begin{equation} \label{eq5_05}
4 \cos^3 \theta-3\cos \theta = \dfrac{3q}{2p} \sqrt{-\dfrac{3}{p}}
\end{equation}\ \\
\ \\
Note that:\\
$\cos 3 \theta =\cos \theta \cos 2\theta -\sin \theta \sin 2\theta=\cos \theta (2\cos^2 \theta-1)-2\cos \theta \sin^2 \theta$\\
\ \\
$\Longrightarrow \hspace{0.5cm} 	\cos 3 \theta = 4 \cos^3 \theta-3\cos \theta)$\\
\ \\
On the basis of this last equality, \eqref {eq5_05} becomes:\\
\begin{equation} \label{eq5_06}
	\cos 3 \theta = \dfrac{3q}{2p} \sqrt{-\dfrac{3}{p}}
\end{equation}\ \\
\ \\
In order for the three roots to be real, it must be:\\
\begin{equation*}
\lvert \dfrac{3q}{2p} \sqrt{-\dfrac{3}{p}} \rvert \leq 1
\end{equation*}\ \\
\ \\
\begin{equation*}
\Longrightarrow \qquad -\dfrac{27q^2}{4p^3}   \leq 1
\end{equation*}\ \\
\begin{equation*}
\Longrightarrow \qquad 27q^2 \leq - 4p^3
\end{equation*}\ \\
\begin{equation} \label{eq5_07}
(\dfrac{q}{2})^2 + ( \dfrac{p}{3})^3 \leq 0
\end{equation}\ \\
\ \\
The above condition \eqref {eq5_07} is more restrictive than \eqref {eq5_04} and is know to be necessary and sufficient for the reality of the three roots from Cardano’s formula.\\
\ \\
From \eqref {eq5_06} you get:\\
\begin{equation*}
\theta = \dfrac{1}{3} \arccos(\dfrac{3q}{2p} \sqrt{-\dfrac{3}{p}})
\end{equation*}\ \\
\ \\
Recalling that:
\begin{equation*}
R=2 \sqrt{-\dfrac{p}{3}}
\end{equation*}\ \\
you get:\\
\ \\
\begin{equation*}
x_{1.2.3}=2\sqrt{-\dfrac{p}{3}} \cos [\dfrac{1}{3} \arccos(\dfrac{3q}{2p} \sqrt{-\dfrac{3}{p}})+\dfrac{2k\pi}{3}] \quad \land \quad k=-1,0,1
\end{equation*}\ \\
\ \\
\end{flushleft}
\centering
q.e.d.\\

\end{document}